\documentclass[dvipsnames,11pt,reqno]{amsart}
\newcommand\ProjectRoot{.}
  \usepackage{amssymb,mathtools,amsthm}
  \usepackage[a4paper,left=3cm,right=3cm,top=3cm,bottom=3cm]{geometry}
  \usepackage[T1]{fontenc}
  \usepackage[utf8]{inputenc}
\usepackage{microtype}
\usepackage{enumitem}

\newcounter{n}
\numberwithin{n}{section}
\theoremstyle{plain}
  \newtheorem{lemma}[n]{Lemma}
  \newtheorem{theorem}[n]{Theorem}
  
  \newtheorem{corollary}[n]{Corollary}
\theoremstyle{definition}
  \newtheorem{definition}[n]{Definition}

  \usepackage{tikz}
  \usepackage{subcaption}

\usepackage{hyperref}
\definecolor{colorlinks}{RGB}{0, 24, 168}
\definecolor{colorcites}{RGB}{124, 10, 2}
\hypersetup{
    colorlinks=true,
    linkcolor=colorlinks,
    citecolor=colorcites,
    urlcolor=colorlinks,
    pdfborder={0 0 0}
}

\renewcommand\phi\varphi
\renewcommand\epsilon\varepsilon
\DeclareMathSymbol{\shortminus}{\mathbin}{AMSa}{"39}
\newcommand\enumindent{}
\begin{document}

 \title[Delocalisation on cubic planar graphs]{Height function delocalisation\\on cubic planar graphs}

 \makeatletter
 \@namedef{subjclassname@2020}{\textup{2020} Mathematics Subject Classification}
 \makeatother

 \author{Piet Lammers}
\address{Institut des Hautes \'Etudes Scientifiques}
\email{lammers@ihes.fr}

\begin{abstract}
The interest is in models of integer-valued height functions on shift-invariant planar graphs whose maximum degree is three.
We prove delocalisation for models induced by convex nearest-neighbour potentials,
under the condition that each potential function is an \emph{excited potential},
that is,
a convex symmetric potential function $V$ with the property that $V(\pm1)\leq V(0)+\log2$.
Examples of such models include the discrete Gaussian and solid-on-solid models at inverse temperature $\beta\leq\log2$,
as well as the uniformly random $K$-Lipschitz function for fixed $K\in\mathbb N$.
In fact, $\beta V$ is an excited potential for any convex symmetric potential function $V$ whenever $\beta$ is sufficiently small.
To arrive at the result, we develop a new technique for symmetry breaking,
and then study the geometric percolation properties of sets of the form $\{\phi\geq a\}$ and $\{\phi\leq a\}$,
where $\phi$ is the random height function and $a$ a constant.
Along the same lines, we derive delocalisation for models induced by convex symmetric nearest-neighbour potentials which force the parity of the height of neighbouring vertices to be distinct.
This includes models of uniformly random graph homomorphisms on the honeycomb lattice and the truncated square tiling,
as well as on the same graphs with each edge replaced by $N$ edges linked in series.
The latter resembles cable-graph constructions which appear in the analysis of the Gaussian free field.

\end{abstract}

\maketitle
\section{Introduction}

Height functions have an important role within statistical physics:
they appear naturally in the study of several lattice models,
such as the dimer models,
percolation models,
the Ising model,
and the six-vertex model,
and height functions are also increasingly studied in their own right.
A natural class of height function models is the class generated by \emph{convex nearest-neighbour potentials},
and we restrict our attention to this class,
even though results have been obtained for non-convex potentials
as well as for non-pair-interactions.
Sheffield was the first to study this class in its full generality in his seminal work \emph{Random Surfaces}~\cite{SHEFFIELD}, where many fundamental properties are derived.
There is a dichotomy for the macroscopic behaviour of each convex nearest-neighbour model:
either the variance of the height difference $\phi(y)-\phi(x)$ remains bounded uniformly over the choice of 
the vertices $x$ and $y$, the so-called \emph{localised} or \emph{smooth} phase,
or the variance of $\phi(y)-\phi(x)$ is unbounded as the distance from $x$ to $y$ grows large, the \emph{delocalised} or \emph{rough} phase.
In dimension two, either class in the dichotomy is nonempty,
and there is not, at present, a general strategy to decide on the class that a model belongs to:
existing results on localisation and delocalisation rely on \emph{ad hoc} arguments.
For a parametrised family of models,
this dichotomy often appears as a form of phase transition:
the six-vertex model with $a=b=1$ and $c\geq 1$, for example,
is localised for $c>2$, and delocalised for $c\leq 2$. 

In dimension two, it is furthermore conjectured that delocalised models have the Gaussian free field as their scaling limit.
This extends the previously mentioned dichotomy;
the scaling limit of a localised model is trivial (or, alternatively, a Gaussian free field of zero variance).
The conjecture implies that a delocalised model must delocalise 
\emph{logarithmically}:
that the variance 
of $\phi(y)-\phi(x)$ grows logarithmically in the distance from $x$ to $y$.
This implied statement has been confirmed for several models,
and no examples to the contrary have been found to date.
Full convergence to the Gaussian free field has been confirmed only for the dimer model, owing to its integrability~\cite{kenyon2001,kenyon2008height,giuliani2017}.

Let us give an overview of delocalisation results for integer-valued height function models
other than the dimer model.
Delocalisation was first proved for the 
discrete Gaussian model and for the solid-on-solid model
at high temperature,
in the landmark article of Fr\"ohlich and Spencer
on the Berezinskii-Kosterlitz-Thouless transition in the XY model, the Villain model,
and the clock model~\cite{frohlich1981}.
The proof relies on a connection with the Coulomb gas.
The six-vertex model with $a=b=1$ and $c\geq 1$
has been shown to localise for $c>2$~\cite{duminil2016discontinuity,GlazmanPeled},
and to delocalise for $c\leq 2$.
Delocalisation was first proven
at several distinct points of the phase diagram~\cite{Dubedat,Giuliani_2017,continuity_rcm_phase_transition,CHAND,GlazmanPeled},
and more recently for the entire range~\cite{in_preparation}.
Delocalisation of the Loop~$O(2)$ model,
that is, the random one-Lipschitz function 
on the vertices of the triangular lattice,
was first established for the parameter $x=1/\sqrt 2$ in~\cite{duminil2017macroscopic},
and later for the parameter $x=1$, the uniform case, in~\cite{glazman2018uniform}.
It is known that the delocalisation is logarithmic
for all examples that were mentioned thus far.
For the six-vertex model and the random Lipschitz function on the triangular lattice,
the proof of this fact relies on Russo-Seymour-Welsh-type
arguments.
(See~\cite{duminil2019logarithmic} for a separate proof of logarithmic delocalisation of square ice.)
Let us finally mention that for planar models,
delocalisation has been established by Sheffield for measures \emph{at a slope},
whenever the slope is not in the dual lattice of the invariance lattice of the underlying model~\cite{SHEFFIELD}.
The proof of this fact relies on \emph{cluster swapping},
which is a versatile method in the analysis of all convex nearest-neighbour models.
We also refer to~\cite{CHAND}, which provides a much more concise proof of this result
for the special case of uniformly random graph homomorphisms from $\mathbb Z^2$ to $\mathbb Z$.

It appears to be the case that the proofs of recent delocalisation results
(for the six-vertex model, the random Lipschitz function on the triangular lattice,
and the general result of Sheffield)
appeal directly to the planar nature of the underlying graph in order to demonstrate that certain random sets of vertices cannot (simultaneously) percolate.
These random sets are often defined in terms of some transformation of the height function $\phi$,
which is chosen to employ the specific duality properties of that specific setting.
We take the exact same route here,
although the choice of random sets is perhaps the most natural one
for the setting of height functions:
we shall directly study the percolative behaviour of the sets
$\{\phi\geq a\}$ and $\{\phi\leq a\}$ for $a$ constant.
This approach is very general,
and it is possible that it leads to similar results for other models when combined with more specialised arguments.
The work of Sheffield~\cite{SHEFFIELD} plays an important role in the proof.

\section{Main results}

We study models of height functions on planar graphs.
In order to prove delocalisation,
we require both the graph as well as the height function model to exhibit certain symmetries.
Let us start with an appropriate definition for the set of graphs that are studied.

\begin{definition}[shift-invariant planar graphs]
	A \emph{shift-invariant planar graph} is a pair $(\mathbb G,\mathcal L)$ where $\mathcal L$ is a lattice in $\mathbb R^2$ and $\mathbb G=(\mathbb V,\mathbb E)$ a connected graph which has a locally finite planar embedding in $\mathbb R^2$ that is invariant under the action of $\mathcal L$.
	The terms \emph{bipartite} and \emph{maximum degree} have their usual meaning;
	in the bipartite case, it is tacitly understood that $\mathcal L$ is chosen such that the two parts of the bipartition are invariant sets.
\end{definition}

\begin{figure}
\centering

\newcommand\customscale{0.8}
\newcommand\customspacing{\hspace{1cm}}

\newcommand\normalscale{\customscale * 1 cm}
\newcommand\correctionscale{\customscale * 0.8135922745 cm}

\begin{subfigure}{.3\textwidth}
\centering
\begin{tikzpicture}[x=\normalscale,y=\normalscale]
\clip (-0.2886751345933333,-0.8333333333333333) rectangle (2.8867513459333334,3.3333333333333335);
\draw (-0.86602540378, -0.5) -- (0.0, 0.0) -- (0.0, 1.0) -- (-0.86602540378, 1.5) -- (0.0, 1.0) -- (0.86602540378, 1.5) -- (1.73205080756, 1.0) -- (1.73205080756, 0.0) -- (0.86602540378, -0.5) -- (0.86602540378, -1.5) -- (0.86602540378, -0.5) -- (0.0, 0.0);
\draw[fill]  (-0.86602540378,-0.5) circle (1pt);
\draw[fill]  (0.0,0.0) circle (1pt);
\draw[fill]  (0.0,1.0) circle (1pt);
\draw[fill]  (-0.86602540378,1.5) circle (1pt);
\draw[fill]  (0.0,1.0) circle (1pt);
\draw[fill]  (0.86602540378,1.5) circle (1pt);
\draw[fill]  (1.73205080756,1.0) circle (1pt);
\draw[fill]  (1.73205080756,0.0) circle (1pt);
\draw[fill]  (0.86602540378,-0.5) circle (1pt);
\draw[fill]  (0.86602540378,-1.5) circle (1pt);
\draw[fill]  (0.86602540378,-0.5) circle (1pt);
\draw[fill]  (0.0,0.0) circle (1pt);
\draw (0.86602540378, 1.5) -- (0.86602540378, 2.5) -- (1.73205080756, 3.0) -- (1.73205080756, 4.0) -- (1.73205080756, 3.0) -- (2.59807621134, 2.5) -- (3.46410161512, 3.0) -- (2.59807621134, 2.5) -- (2.59807621134, 1.5) -- (3.46410161512, 1.0) -- (2.59807621134, 1.5) -- (1.7320508075599998, 1.0);
\draw[fill]  (0.86602540378,1.5) circle (1pt);
\draw[fill]  (0.86602540378,2.5) circle (1pt);
\draw[fill]  (1.73205080756,3.0) circle (1pt);
\draw[fill]  (1.73205080756,4.0) circle (1pt);
\draw[fill]  (1.73205080756,3.0) circle (1pt);
\draw[fill]  (2.59807621134,2.5) circle (1pt);
\draw[fill]  (3.46410161512,3.0) circle (1pt);
\draw[fill]  (2.59807621134,2.5) circle (1pt);
\draw[fill]  (2.59807621134,1.5) circle (1pt);
\draw[fill]  (3.46410161512,1.0) circle (1pt);
\draw[fill]  (2.59807621134,1.5) circle (1pt);
\draw[fill]  (1.7320508075599998,1.0) circle (1pt);
\draw (0.86602540378, 2.5) -- (0.0, 3.0) -- (-0.86602540378, 2.5) -- (0.0, 3.0) -- (0.0, 4.0);
\draw[fill]  (0.86602540378,2.5) circle (1pt);
\draw[fill]  (0.0,3.0) circle (1pt);
\draw[fill]  (-0.86602540378,2.5) circle (1pt);
\draw[fill]  (0.0,3.0) circle (1pt);
\draw[fill]  (0.0,4.0) circle (1pt);
\draw (1.73205080756, 0.0) -- (2.59807621134, -0.5) -- (3.46410161512, 0.0) -- (2.59807621134, -0.5) -- (2.59807621134, -1.5);
\draw[fill]  (1.73205080756,0.0) circle (1pt);
\draw[fill]  (2.59807621134,-0.5) circle (1pt);
\draw[fill]  (3.46410161512,0.0) circle (1pt);
\draw[fill]  (2.59807621134,-0.5) circle (1pt);
\draw[fill]  (2.59807621134,-1.5) circle (1pt);
\end{tikzpicture}
\caption{}
\end{subfigure}
\begin{subfigure}{.3\textwidth}
\centering
\begin{tikzpicture}[x=\correctionscale,y=\correctionscale]
\clip (-0.5,-1.20710678119) rectangle (4.62132034357,3.91421356238);
\draw (0, 0) -- (-1.0, 0.0) -- (0.0, 0.0) -- (0.0, 1.0) -- (-1.0, 1.0) -- (0.0, 1.0) -- (0.70710678119, 1.70710678119) -- (1.70710678119, 1.70710678119) -- (2.41421356238, 1.0) -- (2.41421356238, 0.0) -- (1.70710678119, -0.70710678119) -- (1.70710678119, -1.70710678119) -- (1.70710678119, -0.70710678119) -- (0.70710678119, -0.70710678119) -- (0.70710678119, -1.70710678119) -- (0.70710678119, -0.70710678119) -- (0.0, 0.0);
\draw[fill]  (0,0) circle (1pt);
\draw[fill]  (-1.0,0.0) circle (1pt);
\draw[fill]  (0.0,0.0) circle (1pt);
\draw[fill]  (0.0,1.0) circle (1pt);
\draw[fill]  (-1.0,1.0) circle (1pt);
\draw[fill]  (0.0,1.0) circle (1pt);
\draw[fill]  (0.70710678119,1.70710678119) circle (1pt);
\draw[fill]  (1.70710678119,1.70710678119) circle (1pt);
\draw[fill]  (2.41421356238,1.0) circle (1pt);
\draw[fill]  (2.41421356238,0.0) circle (1pt);
\draw[fill]  (1.70710678119,-0.70710678119) circle (1pt);
\draw[fill]  (1.70710678119,-1.70710678119) circle (1pt);
\draw[fill]  (1.70710678119,-0.70710678119) circle (1pt);
\draw[fill]  (0.70710678119,-0.70710678119) circle (1pt);
\draw[fill]  (0.70710678119,-1.70710678119) circle (1pt);
\draw[fill]  (0.70710678119,-0.70710678119) circle (1pt);
\draw[fill]  (0.0,0.0) circle (1pt);
\draw (1.70710678119, 1.70710678119) -- (1.70710678119, 2.7071067811900003) -- (2.41421356238, 3.4142135623800005) -- (2.41421356238, 4.4142135623800005) -- (2.41421356238, 3.4142135623800005) -- (3.41421356238, 3.4142135623800005) -- (3.41421356238, 4.4142135623800005) -- (3.41421356238, 3.4142135623800005) -- (4.12132034357, 2.7071067811900003) -- (5.12132034357, 2.7071067811900003) -- (4.12132034357, 2.7071067811900003) -- (4.12132034357, 1.7071067811900003) -- (5.12132034357, 1.7071067811900003) -- (4.12132034357, 1.7071067811900003) -- (3.4142135623799996, 1.0000000000000002) -- (2.4142135623799996, 1.0000000000000002);
\draw[fill]  (1.70710678119,1.70710678119) circle (1pt);
\draw[fill]  (1.70710678119,2.7071067811900003) circle (1pt);
\draw[fill]  (2.41421356238,3.4142135623800005) circle (1pt);
\draw[fill]  (2.41421356238,4.4142135623800005) circle (1pt);
\draw[fill]  (2.41421356238,3.4142135623800005) circle (1pt);
\draw[fill]  (3.41421356238,3.4142135623800005) circle (1pt);
\draw[fill]  (3.41421356238,4.4142135623800005) circle (1pt);
\draw[fill]  (3.41421356238,3.4142135623800005) circle (1pt);
\draw[fill]  (4.12132034357,2.7071067811900003) circle (1pt);
\draw[fill]  (5.12132034357,2.7071067811900003) circle (1pt);
\draw[fill]  (4.12132034357,2.7071067811900003) circle (1pt);
\draw[fill]  (4.12132034357,1.7071067811900003) circle (1pt);
\draw[fill]  (5.12132034357,1.7071067811900003) circle (1pt);
\draw[fill]  (4.12132034357,1.7071067811900003) circle (1pt);
\draw[fill]  (3.4142135623799996,1.0000000000000002) circle (1pt);
\draw[fill]  (2.4142135623799996,1.0000000000000002) circle (1pt);
\draw (0.70710678119, 1.70710678119) -- (0.70710678119, 2.7071067811900003) -- (1.70710678119, 2.7071067811900003) -- (0.70710678119, 2.7071067811900003) -- (0.0, 3.4142135623800005) -- (-1.0, 3.4142135623800005) -- (0.0, 3.4142135623800005) -- (0.0, 4.4142135623800005);
\draw[fill]  (0.70710678119,1.70710678119) circle (1pt);
\draw[fill]  (0.70710678119,2.7071067811900003) circle (1pt);
\draw[fill]  (1.70710678119,2.7071067811900003) circle (1pt);
\draw[fill]  (0.70710678119,2.7071067811900003) circle (1pt);
\draw[fill]  (0.0,3.4142135623800005) circle (1pt);
\draw[fill]  (-1.0,3.4142135623800005) circle (1pt);
\draw[fill]  (0.0,3.4142135623800005) circle (1pt);
\draw[fill]  (0.0,4.4142135623800005) circle (1pt);
\draw (2.41421356238, 0.0) -- (3.41421356238, 0.0) -- (3.41421356238, 1.0) -- (3.41421356238, 0.0) -- (4.12132034357, -0.70710678119) -- (5.12132034357, -0.70710678119) -- (4.12132034357, -0.70710678119) -- (4.12132034357, -1.70710678119) -- (4.12132034357, -0.70710678119);
\draw[fill]  (2.41421356238,0.0) circle (1pt);
\draw[fill]  (3.41421356238,0.0) circle (1pt);
\draw[fill]  (3.41421356238,1.0) circle (1pt);
\draw[fill]  (3.41421356238,0.0) circle (1pt);
\draw[fill]  (4.12132034357,-0.70710678119) circle (1pt);
\draw[fill]  (5.12132034357,-0.70710678119) circle (1pt);
\draw[fill]  (4.12132034357,-0.70710678119) circle (1pt);
\draw[fill]  (4.12132034357,-1.70710678119) circle (1pt);
\draw[fill]  (4.12132034357,-0.70710678119) circle (1pt);
\end{tikzpicture}
\caption{}
\end{subfigure}
\begin{subfigure}{.3\textwidth}
\centering
\begin{tikzpicture}[x=\normalscale,y=\normalscale]
\clip (-0.2886751345933333,-0.8333333333333333) rectangle (2.8867513459333334,3.3333333333333335);
\draw (-0.86602540378, -0.5) -- (-0.43301270189, -0.25) -- (0.0, 0.0) -- (0.0, 0.5) -- (0.0, 1.0) -- (-0.43301270189, 1.25) -- (-0.86602540378, 1.5) -- (-0.43301270189, 1.25) -- (0.0, 1.0) -- (0.43301270189, 1.25) -- (0.86602540378, 1.5) -- (1.29903810567, 1.25) -- (1.73205080756, 1.0) -- (1.73205080756, 0.5) -- (1.73205080756, 0.0) -- (1.29903810567, -0.25) -- (0.8660254037799999, -0.5) -- (0.8660254037799999, -1.0) -- (0.8660254037799999, -1.5) -- (0.8660254037799999, -1.0) -- (0.8660254037799999, -0.5) -- (0.4330127018899999, -0.25) -- (-1.1102230246251565e-16, 0.0);
\draw[fill]  (-0.86602540378,-0.5) circle (1pt);
\draw[fill]  (-0.43301270189,-0.25) circle (1pt);
\draw[fill]  (0.0,0.0) circle (1pt);
\draw[fill]  (0.0,0.5) circle (1pt);
\draw[fill]  (0.0,1.0) circle (1pt);
\draw[fill]  (-0.43301270189,1.25) circle (1pt);
\draw[fill]  (-0.86602540378,1.5) circle (1pt);
\draw[fill]  (-0.43301270189,1.25) circle (1pt);
\draw[fill]  (0.0,1.0) circle (1pt);
\draw[fill]  (0.43301270189,1.25) circle (1pt);
\draw[fill]  (0.86602540378,1.5) circle (1pt);
\draw[fill]  (1.29903810567,1.25) circle (1pt);
\draw[fill]  (1.73205080756,1.0) circle (1pt);
\draw[fill]  (1.73205080756,0.5) circle (1pt);
\draw[fill]  (1.73205080756,0.0) circle (1pt);
\draw[fill]  (1.29903810567,-0.25) circle (1pt);
\draw[fill]  (0.8660254037799999,-0.5) circle (1pt);
\draw[fill]  (0.8660254037799999,-1.0) circle (1pt);
\draw[fill]  (0.8660254037799999,-1.5) circle (1pt);
\draw[fill]  (0.8660254037799999,-1.0) circle (1pt);
\draw[fill]  (0.8660254037799999,-0.5) circle (1pt);
\draw[fill]  (0.4330127018899999,-0.25) circle (1pt);
\draw[fill]  (-1.1102230246251565e-16,0.0) circle (1pt);
\draw (0.86602540378, 1.5) -- (0.86602540378, 2.0) -- (0.86602540378, 2.5) -- (1.29903810567, 2.75) -- (1.73205080756, 3.0) -- (1.73205080756, 3.5) -- (1.73205080756, 4.0) -- (1.73205080756, 3.5) -- (1.73205080756, 3.0) -- (2.16506350945, 2.75) -- (2.59807621134, 2.5) -- (3.03108891323, 2.75) -- (3.46410161512, 3.0) -- (3.03108891323, 2.75) -- (2.59807621134, 2.5) -- (2.59807621134, 2.0) -- (2.59807621134, 1.5) -- (3.03108891323, 1.25) -- (3.46410161512, 1.0) -- (3.03108891323, 1.25) -- (2.59807621134, 1.5) -- (2.16506350945, 1.25) -- (1.7320508075599998, 1.0);
\draw[fill]  (0.86602540378,1.5) circle (1pt);
\draw[fill]  (0.86602540378,2.0) circle (1pt);
\draw[fill]  (0.86602540378,2.5) circle (1pt);
\draw[fill]  (1.29903810567,2.75) circle (1pt);
\draw[fill]  (1.73205080756,3.0) circle (1pt);
\draw[fill]  (1.73205080756,3.5) circle (1pt);
\draw[fill]  (1.73205080756,4.0) circle (1pt);
\draw[fill]  (1.73205080756,3.5) circle (1pt);
\draw[fill]  (1.73205080756,3.0) circle (1pt);
\draw[fill]  (2.16506350945,2.75) circle (1pt);
\draw[fill]  (2.59807621134,2.5) circle (1pt);
\draw[fill]  (3.03108891323,2.75) circle (1pt);
\draw[fill]  (3.46410161512,3.0) circle (1pt);
\draw[fill]  (3.03108891323,2.75) circle (1pt);
\draw[fill]  (2.59807621134,2.5) circle (1pt);
\draw[fill]  (2.59807621134,2.0) circle (1pt);
\draw[fill]  (2.59807621134,1.5) circle (1pt);
\draw[fill]  (3.03108891323,1.25) circle (1pt);
\draw[fill]  (3.46410161512,1.0) circle (1pt);
\draw[fill]  (3.03108891323,1.25) circle (1pt);
\draw[fill]  (2.59807621134,1.5) circle (1pt);
\draw[fill]  (2.16506350945,1.25) circle (1pt);
\draw[fill]  (1.7320508075599998,1.0) circle (1pt);
\draw (0.86602540378, 2.5) -- (0.43301270189, 2.75) -- (0.0, 3.0) -- (-0.43301270189, 2.75) -- (-0.86602540378, 2.5) -- (-0.43301270189, 2.75) -- (0.0, 3.0) -- (0.0, 3.5) -- (0.0, 4.0);
\draw[fill]  (0.86602540378,2.5) circle (1pt);
\draw[fill]  (0.43301270189,2.75) circle (1pt);
\draw[fill]  (0.0,3.0) circle (1pt);
\draw[fill]  (-0.43301270189,2.75) circle (1pt);
\draw[fill]  (-0.86602540378,2.5) circle (1pt);
\draw[fill]  (-0.43301270189,2.75) circle (1pt);
\draw[fill]  (0.0,3.0) circle (1pt);
\draw[fill]  (0.0,3.5) circle (1pt);
\draw[fill]  (0.0,4.0) circle (1pt);
\draw (1.73205080756, 0.0) -- (2.16506350945, -0.25) -- (2.59807621134, -0.5) -- (3.03108891323, -0.25) -- (3.46410161512, 0.0) -- (3.03108891323, -0.25) -- (2.59807621134, -0.5) -- (2.59807621134, -1.0) -- (2.59807621134, -1.5);
\draw[fill]  (1.73205080756,0.0) circle (1pt);
\draw[fill]  (2.16506350945,-0.25) circle (1pt);
\draw[fill]  (2.59807621134,-0.5) circle (1pt);
\draw[fill]  (3.03108891323,-0.25) circle (1pt);
\draw[fill]  (3.46410161512,0.0) circle (1pt);
\draw[fill]  (3.03108891323,-0.25) circle (1pt);
\draw[fill]  (2.59807621134,-0.5) circle (1pt);
\draw[fill]  (2.59807621134,-1.0) circle (1pt);
\draw[fill]  (2.59807621134,-1.5) circle (1pt);
\end{tikzpicture}
\caption{}
\end{subfigure}

\caption{%
Examples of suitable graphs:
(a) the honeycomb lattice,
(b) the truncated square tiling,
and (c) the honeycomb lattice with edges replaced by series of edges.
}
\label{fig:graphs}

\end{figure}
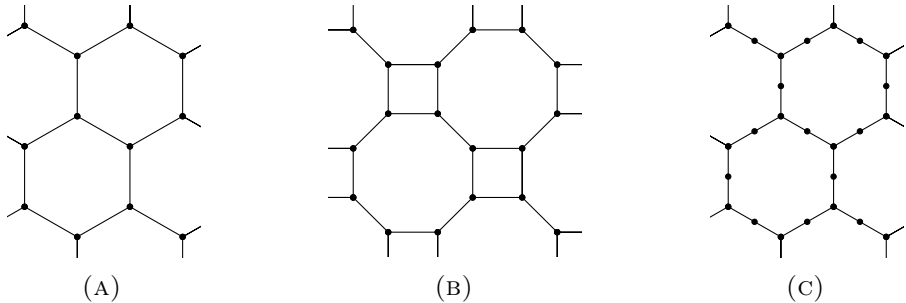

Of particular interest are the shift-invariant planar graphs of maximum degree three.
Prime examples of such graphs include the honeycomb lattice
and the truncated square tiling (sometimes also called the \emph{octagonal tiling} or \emph{Mediterranean tiling}), see Figure~\ref{fig:graphs}.
The pair $(\mathbb G,\mathcal L)$ shall denote a shift-invariant planar graph throughout this article.
The following definition introduces the standard terminology of height functions.

\begin{definition}[height functions]
	A \emph{height function} is an integer-valued function on the vertex set $\mathbb V$.
	Write $\Omega=\mathbb Z^\mathbb V$ for the set of height functions,
	and write $\mathcal F$ for the Borel $\sigma$-algebra of the product topology on $\Omega$.
	Identify each element $\theta\in\mathcal L$ with the map $\mathbb V\to\mathbb V,\,x\mapsto x+\theta$.
	For $A\in\mathcal F$ and $\theta\in\mathcal L$,
	let $\theta A$ denote the event $\{\phi\circ\theta:\phi\in A\}\in\mathcal F$.
	Write $\mathcal P(\Omega,\mathcal F)$ for the set of probability measures on $(\Omega,\mathcal F)$, and write $\mathcal P_\mathcal L(\Omega,\mathcal F)$ for the set of measures which are furthermore \emph{shift-invariant}, that is,
	measures $\mu$ which satisfy $\mu\circ\theta=\mu$ for any $\theta\in\mathcal L$.
\end{definition}

To characterise the model of interest,
we must introduce potential functions which encode how a random height function $\phi$ behaves on a finite set $\Lambda\subset\mathbb V$, conditional on the values of $\phi$ on the complement of $\Lambda$.
We shall first introduce the standard notation in the following definition,
before describing exactly the class of potentials that are analysed here.

\begin{definition}[potential functions, Hamiltonians, Gibbs measures]
\label{definition:Hamiltonian}
Write $\lambda$ for the counting measure on $\mathbb Z$.
A \emph{potential function} is a function $V:\mathbb Z\to\mathbb R\cup\{\infty\}$ with the property that the \emph{edge transition distribution} given by $e^{-V}\lambda$ is a nontrivial finite measure.
The potential function $V$ is called \emph{convex} whenever it is convex as a function over $\mathbb Z$,
that is, if $V((1-t)a+tb)\leq (1-t)V(a)+tV(b)$
for any integers $a$ and $b$ and for any $t\in[0,1]$ such that 
$(1-t)a+tb$ is also an integer.
We shall also impose that $V$ is \emph{symmetric}, in the sense that $V(-x)=V(x)$ for any $x$.
For any finite set $\Lambda\subset\mathbb V$,
introduce the associated \emph{Hamiltonian} defined by
\[
	H_\Lambda:
	\Omega\to\mathbb R\cup\{\infty\},\,
	\phi\mapsto\sum\nolimits_{xy\in\mathbb E(\Lambda)}V(\phi(y)-\phi(x)),
\]
where $\mathbb E(\Lambda)$ denotes the set of undirected edges in $\mathbb E$ that have at least one endpoint in $\Lambda$.
A height function $\phi$ is called \emph{admissible} if $H_\Lambda(\phi)$ is finite for any $\Lambda$.
If $\phi$ is admissible and $\Lambda\subset\mathbb V$ finite, then write 
$\gamma_\Lambda(\cdot,\phi)$ for 
the probability measure 
defined by 
\begin{equation}
\label{eq:spec_def}
	\gamma_\Lambda(\cdot,\phi)
	:=
	\frac1{Z_\Lambda(\phi)}e^{-H_\Lambda}(\delta_{\phi|_{\mathbb V\smallsetminus\Lambda}}\times\lambda^\Lambda),
\end{equation}
where $Z_\Lambda(\phi)$ denotes a suitable normalisation constant.
Thus, to sample from $\gamma_\Lambda(\cdot,\phi)$, set first the random height function equal to $\phi$ on the complement of $\Lambda$,
then sample its values on $\Lambda$ proportional to $e^{-H_\Lambda}$.
The family $(\gamma_\Lambda)_\Lambda$ with $\Lambda$ ranging over the finite subsets of $\mathbb V$ is a specification.
A measure $\mu\in\mathcal P(\Omega,\mathcal F)$ is called a \emph{Gibbs measure} if it is supported on admissible configurations,
and if $\mu=\mu\gamma_\Lambda$ for any finite $\Lambda\subset\mathbb V$.
Write $\mathcal G$ for the collection of Gibbs measures,
and $\mathcal G_\mathcal L$ for the collection 
of shift-invariant Gibbs measures,
that is, $\mathcal G_\mathcal L:=\mathcal G\cap\mathcal P_\mathcal L(\Omega,\mathcal F)$.
\end{definition}

Let us now describe the class of potential functions that are studied.

\begin{definition}[excited potentials]
	An \emph{excited potential} is a convex symmetric potential function $V$ with the property that $V(\pm1)\leq V(0)+\log2$.
\end{definition}

The class of excited potentials is large,
and includes
the discrete Gaussian and solid-on-solid models at inverse temperature $\beta\leq\log2$,
as well as the uniformly random $K$-Lipschitz function for fixed $K\in\mathbb N$.
In fact, if $V$ is any convex symmetric potential function with $V(\pm1)<\infty$,
then $\beta V$ is an excited potential for $\beta\leq \log2/(V(\pm 1)-V(0))$.
With the definition of excited potentials in place,
we are ready to state the main result.

\begin{theorem}[delocalisation]
	\label{theorem:delocalisation_even}
	Let $(\mathbb G,\mathcal L)$ be a shift-invariant planar graph
	of maximum degree three, and
	let $V$ denote an excited potential.
	Then the associated height function model delocalises,
	in the sense that the set $\mathcal G_\mathcal L$ is empty.
\end{theorem}

The fundamental feature of excited potentials is that they facilitate a form of symmetry breaking on the edges on which $\phi$ is constant.
In fact, the case for delocalisation is much simpler if the potential directly prohibits such edges from appearing.
This motivates the following definition.

\begin{definition}[parity potentials]
	A \emph{parity potential} is a symmetric potential function $V$ which satisfies $V(x)=\infty$ for any even integer $x$,
	and whose restriction to the odd integers is convex.
\end{definition}

Indeed, the definition implies that any admissible height function $\phi$ has the property that $\phi(x)$ and $\phi(y)$ have a different parity for any edge $xy\in\mathbb E$,
which implies in particular that the two values cannot be equal.
The graph $\mathbb G$ must thus be bipartite when working with parity potentials, or no admissible height functions would exist.
The prime example of a parity potential is the potential function 
defined by $V(x):=\infty\cdot 1_{|x|\neq 1}$,
which induces a model of uniformly random graph homomorphisms.

\begin{theorem}[delocalisation]
	\label{thm:delocalisation_main}
	Let $(\mathbb G,\mathcal L)$ denote a bipartite shift-invariant planar graph
	of maximum degree three.
	Let $V$ be a parity potential.
	Then the associated height function model delocalises,
	in the sense that the set $\mathcal G_\mathcal L$ is empty.
\end{theorem}

The theorem thus includes models of uniformly random graph homomorphisms 
on the honeycomb lattice and the truncated square tiling.
Another interesting construction is the following:
fix a natural number $N$,
and replace each edge of the honeycomb lattice
by $N$ edges which are linked in series;
see Figure~\ref{fig:graphs}.
The previous theorem applies also to the model of
uniformly random graph homomorphisms on this expanded graph.
This model is reminiscent of cable-graph constructions which appear in the analysis of the (discrete) Gaussian free field,
see for example the work of Lupu~\cite{LUPU}.
The same construction works also for other graphs; it is not restricted to the honeycomb lattice. 
Remark that Theorem~\ref{thm:delocalisation_main} is included in~\cite{CHAND}
whenever there exists an automorphism of $\mathbb G$ that interchanges the two parts of the bipartition of $\mathbb G$.

The proofs of Theorems~\ref{theorem:delocalisation_even} and~\ref{thm:delocalisation_main}
do not rely on the fact that the potential function $V$ is the same for each edge;
the theorems remain valid when replacing the Hamiltonian by
\[
	H_\Lambda:
	\Omega\to\mathbb R\cup\{\infty\},\,
	\phi\mapsto\sum\nolimits_{xy\in\mathbb E(\Lambda)}V_{xy}(\phi(y)-\phi(x)),
\]
where $xy\mapsto V_{xy}$ is a shift-invariant assignment of 
excited potentials or parity potentials to the edges of the graph $\mathbb G$.
It is even possible to mix the two potential types,
up to a small technical modification.
More precisely, one may mix parity potentials with \emph{even} excited potentials:
potential functions $V$ such that $x\mapsto V(2x)$ is an excited potential,
and which satisfy $V(x)=\infty$ for any odd integer $x$.

\section{Proof strategy}
\label{section:strategy}

The proof runs by contradiction.
Throughout the remainder of this article, $\mathbb G$, $\mathcal L$, and $V$ are fixed,
and---in order to arrive at a contradiction---we shall suppose that $\mathcal G_\mathcal L$ is nonempty, and fix some shift-invariant 
Gibbs measure $\mu\in\mathcal G_\mathcal L$.
It is well-known that $\mu$ may be decomposed into ergodic components which are also shift-invariant Gibbs measures,
and by doing so and choosing one such component to replace $\mu$,
we may assume without loss of generality that $\mu$ is itself ergodic.
This means that $\mu(A)\in\{0,1\}$ for any event 
$A$ which satisfies $\theta A=A$
for all $\theta\in\mathcal L$.

\begin{enumerate}
    \item\label{step1} First, we recall the key result from the work of Sheffield, namely
    that the measure $\mu$ introduced above must necessarily
    be extremal. We shall use two 
    corollaries of that result. The first corollary is the FKG lattice condition which
    is satisfied by the height function $\phi$ in the measure $\mu$.
    The second
    corollary consists of the fact that each height $\phi(x)$ 
    is integrable
    in the measure $\mu$.
    The second corollary plays an important role in the derivation of the
    contradiction:
    in the remainder of the proof, we argue that actually $\mu(\phi(x))$ does not
    equal any real number.
    \item\label{step2} Next, we prove a simple lemma, which for any constant $a\in\mathbb Z$ 
    and for any vertex $x$ relates the percolative behaviour
    of the sets $\{\phi\geq a\}$ and $\{\phi\leq a\}$ to the value of $\mu(\phi(x))$.
    If $G$ is any graph and $A$ a random subset of its vertices, then write 
    $X_G(A)$ for the event that the complement of $A$ does not contain an infinite
    cluster.
    On this event, we think of $A$ as preventing its complement from percolating.
    Note that the events $X_\mathbb G(\{\phi\geq a\})$ and $X_\mathbb G(\{\phi\leq a\})$
    are shift-invariant and therefore satisfy a zero-one law in the ergodic measure $\mu$.
    The lemma asserts that $\mu(\phi(x))\geq a$ whenever the event $X_\mathbb G(\{\phi\geq a\})$
    occurs almost surely,
    and that $\mu(\phi(x))\leq a$ whenever the event $X_\mathbb G(\{\phi\leq a\})$
    occurs almost surely.
    The proof uses a natural exploration process,
    as well as the symmetry of the potential function $V$.
    \item\label{step3} In the following step, we compare the events $X_\mathbb G(\{\phi\leq a-1\})$
    and $X_\mathbb G(\{\phi\geq a\})$.
    If both events occur, then the previous lemma implies
    that $\mu(\phi(x))\leq a-1$ and simultaneously $\mu(\phi(x))\geq a$, a contradiction.
    If both events do not occur, then both $\{\phi\geq a\}$
    and $\{\phi\leq a-1\}$ percolate almost surely.
    Phase coexistence for planar site percolation was, however, ruled out 
    by an argument of Sheffield.
    Thus, we conclude that for any integer $a$, exactly one of $X_\mathbb G(\{\phi\leq a-1\})$
    and $X_\mathbb G(\{\phi\geq a\})$ occurs almost surely,
    and that the other event occurs almost never.
    This also implies that $\mu(\phi(x))$ must be an integer,
    and we assume that $\mu(\phi(x))=0$ without loss of generality.
\end{enumerate}
Observe that up to this point, we did not use our two special conditions:
that the maximum degree of the graph is three,
and that the convex symmetric potential function $V$ is either excited
or a parity potential.
The remainder of the proof is easier for parity potentials,
which is why we specialise to these potentials first.
In the next two steps, we use the geometry of the graph $\mathbb G$
to take the final leap towards the contradiction for parity potentials.
In the last step, we generalise the argument to excited potentials.
\begin{enumerate}[resume]
    \item Planarity is of major importance, because of the phase coexistence result 
    introduced above.
    The only reason for requiring the maximum degree of the planar graph $\mathbb G=(\mathbb V,\mathbb E)$ 
    to be three, is that it also implies planarity for certain derived graphs.
    The first derived graph for which this holds true,
    is the \emph{line graph} of $\mathbb G$, that is, the graph which has 
    $\mathbb E$ as its vertex set, and for which two edges are neighbours 
    if and only if they share a vertex in the original graph $\mathbb G$.
    Since the line graph of $\mathbb G$ is planar, the phase coexistence result
    applies also to \emph{edge percolation} on $\mathbb G$.
    The second derived graph is only defined in the case that $\mathbb G$ is bipartite.
    In that case, let $\mathbb V_1$ denote the set of odd vertices,
    and write $\mathbb G_1$ for the graph on $\mathbb V_1$ in which two vertices 
    are neighbours if and only if they are at a graph distance two in the original 
    graph $\mathbb G$.
    It is straightforward to see that the two derived graphs are planar,
    because two adjacent edges (in the case of the line graph) or two adjacent odd vertices
    (in the bipartite case) always belong to the same face of the original graph $\mathbb G$;
    see Figure~\ref{fig:graphs_expanded}.
    It is also easy to see that the condition on the maximum degree is necessary:
    the derived graphs would clearly not be planar for the square lattice graph.
    \begin{figure}
\centering

\newcommand\customscale{0.8}
\newcommand\customspacing{\hspace{1cm}}

\newcommand\normalscale{\customscale * 1 cm}
\newcommand\correctionscale{\customscale * 0.8135922745 cm}

\begin{subfigure}{.3\textwidth}
\centering
\begin{tikzpicture}[x=\correctionscale,y=\correctionscale]
\clip (-0.5,-1.20710678119) rectangle (4.62132034357,3.91421356238);
\draw (0, 0) -- (-1.0, 0.0) -- (0.0, 0.0) -- (0.0, 1.0) -- (-1.0, 1.0) -- (0.0, 1.0) -- (0.70710678119, 1.70710678119) -- (1.70710678119, 1.70710678119) -- (2.41421356238, 1.0) -- (2.41421356238, 0.0) -- (1.70710678119, -0.70710678119) -- (1.70710678119, -1.70710678119) -- (1.70710678119, -0.70710678119) -- (0.70710678119, -0.70710678119) -- (0.70710678119, -1.70710678119) -- (0.70710678119, -0.70710678119) -- (0.0, 0.0);
\draw[fill]  (0,0) circle (1pt);
\draw[fill]  (-1.0,0.0) circle (1pt);
\draw[fill]  (0.0,0.0) circle (1pt);
\draw[fill]  (0.0,1.0) circle (1pt);
\draw[fill]  (-1.0,1.0) circle (1pt);
\draw[fill]  (0.0,1.0) circle (1pt);
\draw[fill]  (0.70710678119,1.70710678119) circle (1pt);
\draw[fill]  (1.70710678119,1.70710678119) circle (1pt);
\draw[fill]  (2.41421356238,1.0) circle (1pt);
\draw[fill]  (2.41421356238,0.0) circle (1pt);
\draw[fill]  (1.70710678119,-0.70710678119) circle (1pt);
\draw[fill]  (1.70710678119,-1.70710678119) circle (1pt);
\draw[fill]  (1.70710678119,-0.70710678119) circle (1pt);
\draw[fill]  (0.70710678119,-0.70710678119) circle (1pt);
\draw[fill]  (0.70710678119,-1.70710678119) circle (1pt);
\draw[fill]  (0.70710678119,-0.70710678119) circle (1pt);
\draw[fill]  (0.0,0.0) circle (1pt);
\draw (1.70710678119, 1.70710678119) -- (1.70710678119, 2.7071067811900003) -- (2.41421356238, 3.4142135623800005) -- (2.41421356238, 4.4142135623800005) -- (2.41421356238, 3.4142135623800005) -- (3.41421356238, 3.4142135623800005) -- (3.41421356238, 4.4142135623800005) -- (3.41421356238, 3.4142135623800005) -- (4.12132034357, 2.7071067811900003) -- (5.12132034357, 2.7071067811900003) -- (4.12132034357, 2.7071067811900003) -- (4.12132034357, 1.7071067811900003) -- (5.12132034357, 1.7071067811900003) -- (4.12132034357, 1.7071067811900003) -- (3.4142135623799996, 1.0000000000000002) -- (2.4142135623799996, 1.0000000000000002);
\draw[fill]  (1.70710678119,1.70710678119) circle (1pt);
\draw[fill]  (1.70710678119,2.7071067811900003) circle (1pt);
\draw[fill]  (2.41421356238,3.4142135623800005) circle (1pt);
\draw[fill]  (2.41421356238,4.4142135623800005) circle (1pt);
\draw[fill]  (2.41421356238,3.4142135623800005) circle (1pt);
\draw[fill]  (3.41421356238,3.4142135623800005) circle (1pt);
\draw[fill]  (3.41421356238,4.4142135623800005) circle (1pt);
\draw[fill]  (3.41421356238,3.4142135623800005) circle (1pt);
\draw[fill]  (4.12132034357,2.7071067811900003) circle (1pt);
\draw[fill]  (5.12132034357,2.7071067811900003) circle (1pt);
\draw[fill]  (4.12132034357,2.7071067811900003) circle (1pt);
\draw[fill]  (4.12132034357,1.7071067811900003) circle (1pt);
\draw[fill]  (5.12132034357,1.7071067811900003) circle (1pt);
\draw[fill]  (4.12132034357,1.7071067811900003) circle (1pt);
\draw[fill]  (3.4142135623799996,1.0000000000000002) circle (1pt);
\draw[fill]  (2.4142135623799996,1.0000000000000002) circle (1pt);
\draw (0.70710678119, 1.70710678119) -- (0.70710678119, 2.7071067811900003) -- (1.70710678119, 2.7071067811900003) -- (0.70710678119, 2.7071067811900003) -- (0.0, 3.4142135623800005) -- (-1.0, 3.4142135623800005) -- (0.0, 3.4142135623800005) -- (0.0, 4.4142135623800005);
\draw[fill]  (0.70710678119,1.70710678119) circle (1pt);
\draw[fill]  (0.70710678119,2.7071067811900003) circle (1pt);
\draw[fill]  (1.70710678119,2.7071067811900003) circle (1pt);
\draw[fill]  (0.70710678119,2.7071067811900003) circle (1pt);
\draw[fill]  (0.0,3.4142135623800005) circle (1pt);
\draw[fill]  (-1.0,3.4142135623800005) circle (1pt);
\draw[fill]  (0.0,3.4142135623800005) circle (1pt);
\draw[fill]  (0.0,4.4142135623800005) circle (1pt);
\draw (2.41421356238, 0.0) -- (3.41421356238, 0.0) -- (3.41421356238, 1.0) -- (3.41421356238, 0.0) -- (4.12132034357, -0.70710678119) -- (5.12132034357, -0.70710678119) -- (4.12132034357, -0.70710678119) -- (4.12132034357, -1.70710678119) -- (4.12132034357, -0.70710678119);
\draw[fill]  (2.41421356238,0.0) circle (1pt);
\draw[fill]  (3.41421356238,0.0) circle (1pt);
\draw[fill]  (3.41421356238,1.0) circle (1pt);
\draw[fill]  (3.41421356238,0.0) circle (1pt);
\draw[fill]  (4.12132034357,-0.70710678119) circle (1pt);
\draw[fill]  (5.12132034357,-0.70710678119) circle (1pt);
\draw[fill]  (4.12132034357,-0.70710678119) circle (1pt);
\draw[fill]  (4.12132034357,-1.70710678119) circle (1pt);
\draw[fill]  (4.12132034357,-0.70710678119) circle (1pt);
\end{tikzpicture}
\caption{}
\end{subfigure}
\begin{subfigure}{.3\textwidth}
\centering
\begin{tikzpicture}[x=\correctionscale,y=\correctionscale]
\input{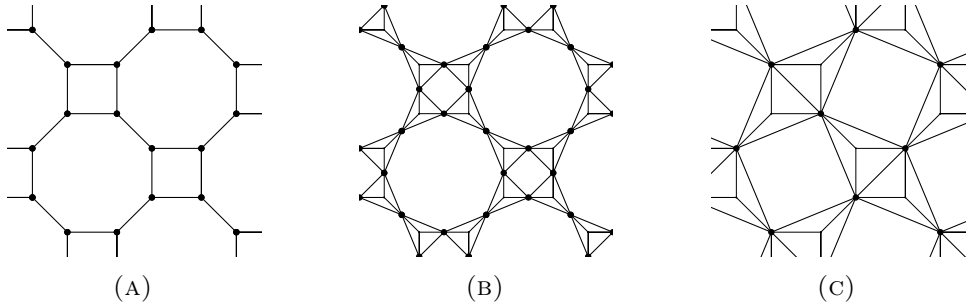}
\end{tikzpicture}
\caption{}
\end{subfigure}
\begin{subfigure}{.3\textwidth}
\centering
\begin{tikzpicture}[x=\correctionscale,y=\correctionscale]
\clip (-0.5,-1.20710678119) rectangle (4.62132034357,3.91421356238);
\draw[very thin] (0, 0) -- (-1.0, 0.0) -- (0.0, 0.0) -- (0.0, 1.0) -- (-1.0, 1.0) -- (0.0, 1.0) -- (0.70710678119, 1.70710678119) -- (1.70710678119, 1.70710678119) -- (2.41421356238, 1.0) -- (2.41421356238, 0.0) -- (1.70710678119, -0.70710678119) -- (1.70710678119, -1.70710678119) -- (1.70710678119, -0.70710678119) -- (0.70710678119, -0.70710678119) -- (0.70710678119, -1.70710678119) -- (0.70710678119, -0.70710678119) -- (0.0, 0.0);
\draw[fill]  (0,0) circle (0.1pt);
\draw[fill]  (-1.0,0.0) circle (0.1pt);
\draw[fill]  (0.0,0.0) circle (0.1pt);
\draw[fill]  (0.0,1.0) circle (0.1pt);
\draw[fill]  (-1.0,1.0) circle (0.1pt);
\draw[fill]  (0.0,1.0) circle (0.1pt);
\draw[fill]  (0.70710678119,1.70710678119) circle (0.1pt);
\draw[fill]  (1.70710678119,1.70710678119) circle (0.1pt);
\draw[fill]  (2.41421356238,1.0) circle (0.1pt);
\draw[fill]  (2.41421356238,0.0) circle (0.1pt);
\draw[fill]  (1.70710678119,-0.70710678119) circle (0.1pt);
\draw[fill]  (1.70710678119,-1.70710678119) circle (0.1pt);
\draw[fill]  (1.70710678119,-0.70710678119) circle (0.1pt);
\draw[fill]  (0.70710678119,-0.70710678119) circle (0.1pt);
\draw[fill]  (0.70710678119,-1.70710678119) circle (0.1pt);
\draw[fill]  (0.70710678119,-0.70710678119) circle (0.1pt);
\draw[fill]  (0.0,0.0) circle (0.1pt);
\draw[very thin] (1.70710678119, 1.70710678119) -- (1.70710678119, 2.7071067811900003) -- (2.41421356238, 3.4142135623800005) -- (2.41421356238, 4.4142135623800005) -- (2.41421356238, 3.4142135623800005) -- (3.41421356238, 3.4142135623800005) -- (3.41421356238, 4.4142135623800005) -- (3.41421356238, 3.4142135623800005) -- (4.12132034357, 2.7071067811900003) -- (5.12132034357, 2.7071067811900003) -- (4.12132034357, 2.7071067811900003) -- (4.12132034357, 1.7071067811900003) -- (5.12132034357, 1.7071067811900003) -- (4.12132034357, 1.7071067811900003) -- (3.4142135623799996, 1.0000000000000002) -- (2.4142135623799996, 1.0000000000000002);
\draw[fill]  (1.70710678119,1.70710678119) circle (0.1pt);
\draw[fill]  (1.70710678119,2.7071067811900003) circle (0.1pt);
\draw[fill]  (2.41421356238,3.4142135623800005) circle (0.1pt);
\draw[fill]  (2.41421356238,4.4142135623800005) circle (0.1pt);
\draw[fill]  (2.41421356238,3.4142135623800005) circle (0.1pt);
\draw[fill]  (3.41421356238,3.4142135623800005) circle (0.1pt);
\draw[fill]  (3.41421356238,4.4142135623800005) circle (0.1pt);
\draw[fill]  (3.41421356238,3.4142135623800005) circle (0.1pt);
\draw[fill]  (4.12132034357,2.7071067811900003) circle (0.1pt);
\draw[fill]  (5.12132034357,2.7071067811900003) circle (0.1pt);
\draw[fill]  (4.12132034357,2.7071067811900003) circle (0.1pt);
\draw[fill]  (4.12132034357,1.7071067811900003) circle (0.1pt);
\draw[fill]  (5.12132034357,1.7071067811900003) circle (0.1pt);
\draw[fill]  (4.12132034357,1.7071067811900003) circle (0.1pt);
\draw[fill]  (3.4142135623799996,1.0000000000000002) circle (0.1pt);
\draw[fill]  (2.4142135623799996,1.0000000000000002) circle (0.1pt);
\draw[very thin] (0.70710678119, 1.70710678119) -- (0.70710678119, 2.7071067811900003) -- (1.70710678119, 2.7071067811900003) -- (0.70710678119, 2.7071067811900003) -- (0.0, 3.4142135623800005) -- (-1.0, 3.4142135623800005) -- (0.0, 3.4142135623800005) -- (0.0, 4.4142135623800005);
\draw[fill]  (0.70710678119,1.70710678119) circle (0.1pt);
\draw[fill]  (0.70710678119,2.7071067811900003) circle (0.1pt);
\draw[fill]  (1.70710678119,2.7071067811900003) circle (0.1pt);
\draw[fill]  (0.70710678119,2.7071067811900003) circle (0.1pt);
\draw[fill]  (0.0,3.4142135623800005) circle (0.1pt);
\draw[fill]  (-1.0,3.4142135623800005) circle (0.1pt);
\draw[fill]  (0.0,3.4142135623800005) circle (0.1pt);
\draw[fill]  (0.0,4.4142135623800005) circle (0.1pt);
\draw[very thin] (2.41421356238, 0.0) -- (3.41421356238, 0.0) -- (3.41421356238, 1.0) -- (3.41421356238, 0.0) -- (4.12132034357, -0.70710678119) -- (5.12132034357, -0.70710678119) -- (4.12132034357, -0.70710678119) -- (4.12132034357, -1.70710678119) -- (4.12132034357, -0.70710678119);
\draw[fill]  (2.41421356238,0.0) circle (0.1pt);
\draw[fill]  (3.41421356238,0.0) circle (0.1pt);
\draw[fill]  (3.41421356238,1.0) circle (0.1pt);
\draw[fill]  (3.41421356238,0.0) circle (0.1pt);
\draw[fill]  (4.12132034357,-0.70710678119) circle (0.1pt);
\draw[fill]  (5.12132034357,-0.70710678119) circle (0.1pt);
\draw[fill]  (4.12132034357,-0.70710678119) circle (0.1pt);
\draw[fill]  (4.12132034357,-1.70710678119) circle (0.1pt);
\draw[fill]  (4.12132034357,-0.70710678119) circle (0.1pt);
\draw (-1.70710678119, -1.70710678119) -- (-1.0, 0.0) -- (0.70710678119, -0.70710678119) -- (0.0, -2.41421356238) -- (-1.70710678119, -1.70710678119);
\draw (0, 1) -- (1.70710678119, 1.70710678119) -- (2.41421356238, 0.0) -- (0.70710678119, -0.70710678119) -- (0.0, 1.0);
\draw (1.70710678119, 1.70710678119) -- (2.41421356238, 3.41421356238) -- (4.12132034357, 2.7071067811900003) -- (3.4142135623799996, 1.0000000000000002) -- (1.7071067811899996, 1.7071067811900003);
\draw (3.41421356238, 4.4142135623800005) -- (5.12132034357, 5.121320343570001) -- (5.82842712476, 3.4142135623800005) -- (4.12132034357, 2.7071067811900003) -- (3.4142135623799996, 4.4142135623800005);
\draw (-1.70710678119, 1.70710678119) -- (-1.0, 3.41421356238) -- (0.70710678119, 2.7071067811900003) -- (0.0, 1.0000000000000002) -- (-1.70710678119, 1.7071067811900003);
\draw (0.0, 4.4142135623800005) -- (1.70710678119, 5.121320343570001) -- (2.41421356238, 3.4142135623800005) -- (0.70710678119, 2.7071067811900003) -- (0.0, 4.4142135623800005);
\draw (1.70710678119, -1.70710678119) -- (2.41421356238, 0.0) -- (4.12132034357, -0.70710678119) -- (3.4142135623799996, -2.41421356238) -- (1.7071067811899996, -1.70710678119);
\draw (3.41421356238, 1.0) -- (5.12132034357, 1.70710678119) -- (5.82842712476, 0.0) -- (4.12132034357, -0.70710678119) -- (3.4142135623799996, 1.0);
\draw (0, 1) -- (-1.0, 0.0);
\draw[fill]  (0,1) circle (1pt);
\draw[fill]  (-1.0,0.0) circle (1pt);
\draw (1.70710678119, 1.70710678119) -- (0.70710678119, 2.7071067811900003);
\draw[fill]  (1.70710678119,1.70710678119) circle (1pt);
\draw[fill]  (0.70710678119,2.7071067811900003) circle (1pt);
\draw (3.41421356238, 4.4142135623800005) -- (2.41421356238, 3.4142135623800005);
\draw[fill]  (3.41421356238,4.4142135623800005) circle (1pt);
\draw[fill]  (2.41421356238,3.4142135623800005) circle (1pt);
\draw (1.70710678119, -1.70710678119) -- (0.70710678119, -0.70710678119);
\draw[fill]  (1.70710678119,-1.70710678119) circle (1pt);
\draw[fill]  (0.70710678119,-0.70710678119) circle (1pt);
\draw (3.41421356238, 1.0) -- (2.41421356238, 0.0);
\draw[fill]  (3.41421356238,1.0) circle (1pt);
\draw[fill]  (2.41421356238,0.0) circle (1pt);
\draw (5.12132034357, 1.70710678119) -- (4.12132034357, 2.7071067811900003);
\draw[fill]  (5.12132034357,1.70710678119) circle (1pt);
\draw[fill]  (4.12132034357,2.7071067811900003) circle (1pt);
\draw (5.12132034357, -1.70710678119) -- (4.12132034357, -0.70710678119);
\draw[fill]  (5.12132034357,-1.70710678119) circle (1pt);
\draw[fill]  (4.12132034357,-0.70710678119) circle (1pt);
\end{tikzpicture}
\caption{}
\end{subfigure}

\caption{%
The truncated square tiling (a) and two derived graphs:
the line graph (b) and the graph of odd vertices (c).
The derived graphs are planar because each edge of the derived graph 
can be drawn on a face of the original graph.
}
\label{fig:graphs_expanded}

\end{figure}

    \item\label{step5} The contradiction for parity potentials is now straightforward.
    Note that the height function $\phi$ takes odd values on the odd vertices almost surely,
    without loss of generality.
    Write $\psi$ for the restriction of $\phi$ to the set of odd vertices.
    As in Step~\ref{step3}, we observe that the phase coexistence result implies that 
    at least one of
    $X_{\mathbb G_1}(\{\psi\geq 1\})$ a $X_{\mathbb G_1}(\{\psi\leq-1\})$
    must occur almost surely
    (say the former, without loss of generality).
    But the event $X_{\mathbb G_1}(\{\psi\geq 1\})$ is included in the event 
    $X_\mathbb G(\{\phi\geq 1\})$, and therefore the latter event must also occur
    almost surely.
    But then $\mu(\phi(x))\geq 1$, which clearly contradicts that $\mu(\phi(x))=0$.
    This completes the proof for parity potentials.

    \item In this final step, we generalise to excited potentials.
    The crucial ingredient in the previous step was that for any edge $xy$
    it almost never occurs that $\phi(x)=\phi(y)$.
    More specifically, for parity potentials it almost never occurs that both $\phi(x)=0$
    and
    $\phi(y)=0$.
    The challenge when generalising to excited potentials is to handle those 
    edges, which now do occur with positive probability.
    Call such edges \emph{zero edges}.

    \enumindent
    To prove delocalisation, we find a way to replace a zero edge by two edges 
    which are linked in series.
    A new vertex appears at the midpoint of the edge,
    and this new vertex will have a height of $\pm\frac12$.
    This means that each (original) edge has at least one vertex with a nonzero 
    height on it.
    The argument is then closed by proceeding more or less as in the previous step:
    by observing that there cannot be phase coexistence for edge percolation 
    and by running the correct exploration process,
    one derives that $|\mu(\phi(x))|\geq\frac12$.
    This contradicts the assumption that $\mu(\phi(x))=0$.
    
    \enumindent
    This new vertex is created as follows.
    Assume that $V(0)=0$ without loss of generality,
    and define $V^*(x):=\log 2\cdot 1_{|x|=1}+\infty\cdot1_{|x|\geq 2}$.
    The relative weight of the height transition at some edge $xy$ is decomposed as
    \[
        e^{-V(h)}
        =e^{-V^*(h)}
            +\left(e^{-V(h)}-e^{-V^*(h)}\right),
        \quad
        \text{where}
        \quad
        h:=\phi(y)-\phi(x).
    \]
    The definition of an excited potential implies that $V^*\geq V$,
    and therefore the summand on the right is always nonnegative.
    Moreover, if $h=0$, then the term on the right does not contribute at all.
    We will argue later that the zero edges (and potentially some others)
    behave as if they have the potential function $V^*$ associated to them rather than 
    the potential function $V$.
    But it is readily observed that the distribution $e^{-V^*}\lambda$ equals
    (up to scaling) the distribution of the sum of two independent fair 
    $\pm\frac12$-valued coin flips.
    This justifies the introduction of the intermediate height (which is essentially 
    the outcome of the first coin flip),
    and leads to the desired contradiction.
\end{enumerate}

\section{Discussion}

\subsection{Interpretation of the conditions in the main result}
Many classical spin lattice models are known to exhibit exponential decay of correlations at a sufficiently high temperature.
For several such models, this can be proven through a coupling with a subcritical percolation model.
The Peierls argument, on the other hand, implies that such models localise at a sufficiently low temperature.
The Peierls argument also applies to height function models:
in the setting of the current article, for example, one may use it to demonstrate
that the model is localised whenever $V(\pm 1)$ is much larger than $V(0)$.
Generic arguments for delocalisation, however, have not been found to date:
there is an intuition that height function models should delocalise at high temperature,
but this intuition has only been verified for a handful of cases.
We argue that this article is a first step towards a generic delocalisation argument.
The first condition---the required bound on $V(\pm 1)-V(0)$---may be interpreted directly as a temperature constraint.
It says that the weight of a height difference of one over an edge should not be too small compared to the weight 
of a constant height.
The second condition---that the maximum degree of $\mathbb G$ is three---is a combinatorial constraint,
but still quite a natural one because a lower maximum degree may be thought of as having a graph with less edges,
which in turn can be thought of as imposing a weaker potential.

\subsection{Improving the conditions in the main result}
The quantitative requirement that $V(\pm1)-V(0)\leq\log2$ was introduced to break symmetry
for the zero edges.
This is followed by a subtle geometrical argument to show that certain sets of sites or bonds
cannot simultaneously percolate.
Therefore it seems that the most natural way to generalise Theorems~\ref{theorem:delocalisation_even} and~\ref{thm:delocalisation_main}
while preserving the overall proof strategy,
would be to improve the percolation argument in order to extend to shift-invariant planar 
graphs which are not cubic (with the same quantitative constraint on $V$).

The most interesting generalisation---in the eyes of the author---would be to extend to the square lattice.
We conjecture that the result holds true on this graph, with the same requirement for $V$.
The reason that the author believes some modification of the proof to work is that one ultimately needs to
rule out phase coexistence for a bond percolation model
(with its complement)
in which the zero edges, whose state is uniformly random and independent of anything else, play an essential role.
This is significant because it is known that uniformly random bond percolation on the square lattice 
does not percolate.
One problem is that the percolation model obtained in this article would 
be self-complementary---and not self-dual---which complicates the route to the result.

An adaptation of the proof will certainly not work on the triangular lattice:
it is known that the Loop~$O(2)$ model is localised for $x=1/\sqrt{3}$,
which corresponds to $V(\pm1)-V(0)=(\log3)/2<\log2$; see~\cite{glazman2018exponential}.
It means that our method for symmetry breaking is fundamentally insufficient
for deriving delocalisation in this case.
This view is consistent with the fact that uniformly random bond percolation on the triangular 
lattice falls into the supercritical regime.
It is perhaps reasonable to expect that convex symmetric height function models on the triangular lattice delocalise 
whenever $V(\pm1)-V(0)\leq(\log 2)/2$,
since this corresponds to the parameter $x=1/\sqrt{2}$ for which delocalisation
was obtained in the Loop~$O(2)$
model in~\cite{duminil2017macroscopic}.

In general, there does not always exist a constant $\beta(\mathbb G)$ such that convex 
symmetric height functions on $\mathbb G$ delocalise whenever $V(\pm1)-V(0)\leq\beta(\mathbb G)$.
As a counterexample, there exist shift-invariant planar graphs on which the uniformly random $1$-Lipschitz
function is localised.
For example, fix $N\in\mathbb N$ large,
start with the hexagonal lattice, and
replace each edge by $N$ edges linked in parallel.
Subsequently replace each edge by two edges linked in series,
as in Figure~\ref{fig:graphs}c.
The $1$-Lipschitz function on this graph is localised due to a Peierls argument.

\subsection{Relation to other works on height functions}
A significant difference between this article and the works of Sheffield~\cite{SHEFFIELD}
and Chandgotia, Peled, Sheffield, and Tassy~\cite{CHAND}
consists in the fact that those works crucially analyse the percolative behaviour 
of the level sets of the random function $\psi-\phi$,
where $\phi$ and $\psi$ denote independent samples from $\mu$,
while the current article (also) considers the level sets of the random function $\phi$
directly.
The latter approach leads to new observations in Steps~\ref{step2} and~\ref{step3} in Section~\ref{section:strategy},
which extend beyond the main result in this article and may be of independent interest.

Theorem~\ref{thm:delocalisation_main} is truly different from the delocalisation result in~\cite{CHAND},
which also involves a parity potential.
In~\cite{CHAND}, the result crucially requires the graph to have an automorphism which interchanges
the two parts of the bipartition of $\mathbb G$.
In the present paper, we never rely on this symmetry.
To see that~\cite{CHAND} genuinely requires the symmetry,
observe that the uniformly random graph homomorphism on the graph described at the end of the previous subsection
is localised for $N$ sufficiently large.

\subsection{Relation to general spin models}
Delocalisation of height function models is known to be related to delocalisation of spin models;
it was already mentioned that the discrete Gaussian model and the solid-on-solid model appear in 
the derivation of the 
Berezinskii-Kosterlitz-Thouless transition.
Since the main result of this article generically includes delocalisation of these two height function models,
one may wonder if delocalisation of other spin lattice models can be derived from the height function delocalisation result.
This problem remains open.

\section{Proof of the main results}

We prove two delocalisation results: one for parity potentials,
and one for excited potentials.
The result for parity potentials follows from a simple geometrical argument.
The result for excited potentials is harder to derive,
and requires the introduction of external randomness to break certain symmetries that appear.
Subsection~\ref{subsec:symmbreaking} describes this technique for symmetry breaking in detail.
Subsection~\ref{subsec:quote_sheffield} states a number of existing results 
from the work of Sheffield,
and corresponds to Step~\ref{step1} of the proof strategy.
Subsection~\ref{subsec:proof_main_theorem} contains the general geometrical argument,
which leads directly to the result for parity potentials.
This corresponds to Steps~\ref{step2}--\ref{step5}.
Subsection~\ref{subsec:finalbigproof} combines the ideas 
of Subsections~\ref{subsec:symmbreaking} and~\ref{subsec:proof_main_theorem},
in order to derive delocalisation for excited potentials.
This is the final step of the proof strategy outlined before.

\subsection{Excited potentials and symmetry breaking}
\label{subsec:symmbreaking}

Essential to the proof of Theorem~\ref{thm:delocalisation_main} is the fact that the two endpoints of an edge cannot have the same height when working with a parity potential.
This implies in particular that the height of at least one of the endpoints of an edge is nonzero.
The same statement is simply false in the context of an excited potential:
it is perfectly possible that two neighbouring vertices have height zero.
In this subsection we introduce a form of symmetry breaking
which allows us to choose a nonzero height for each such edge.

Let $V$ denote an excited potential,
and suppose that $V(0)=0$ without loss of generality
(by adding a constant to the potential if necessary).
Write $V^*$ for the potential function
\[
V^*(x):=
\begin{cases}
0&\text{if $x=0$,}\\
\log 2&\text{if $|x|=1$,}\\
\infty&\text{if $|x|>1$.}
\end{cases}
\]
Then $V\leq V^*$ by definition of an excited potential.
The relative weight $e^{-H_\Lambda}$ of each configuration in the definition 
of $\gamma_\Lambda$ in~\eqref{eq:spec_def} decomposes as a product of the factors
\[
	e^{-V(\phi(y)-\phi(x))}
\]
over the edges $xy\in\mathbb E(\Lambda)$.
We shall further decompose the factor corresponding to a single fixed edge $xy\in\mathbb E(\Lambda)$,
while keeping all other weights the same.
First, decompose the weight corresponding to $xy$ as the sum 
\begin{equation}
\label{eq:weight_decomp}
	e^{-V^*(\phi(y)-\phi(x))}+\left[e^{-V(\phi(y)-\phi(x))}-e^{-V^*(\phi(y)-\phi(x))}\right].
\end{equation}
Note that either term is nonnegative.
Introduce a new random variable,
which indicates if this edge is \emph{excited} or \emph{not excited},
with the weight corresponding to each state on the left and right in~\eqref{eq:weight_decomp}.
Formally, the enriched model can be described as follows:
we study samples $(\phi,\epsilon)\in\mathbb Z^\mathbb V\times\{0,1\}^\mathbb E$,
and the relative weight of each configuration $(\phi,\epsilon)$ is given by multiplying
\[
	1_{\epsilon(xy)=1}e^{-V^*(\phi(y)-\phi(x))}
	+1_{\epsilon(xy)=0}\left[e^{-V(\phi(y)-\phi(x))}-e^{-V^*(\phi(y)-\phi(x))}\right]
\]
over the edges $xy$ involved in the probability kernel.
The value $\epsilon(xy)=1$ indicates that the edge $xy$ is excited.
Let us make two simple and essential observations.
First, the state of $xy$ and the height function $\phi$ 
are independent conditional on $\phi(x)$ and $\phi(y)$.
Second, an edge $xy$ is almost surely excited whenever
$\phi(x)=\phi(y)$.
In particular, this implies that at least one of $\phi(x)$ and $\phi(y)$ is not equal to 
zero almost surely whenever the edge $xy$ is \emph{not} excited.
To complete the picture, we also remark that the edge $xy$ is almost surely not excited whenever 
$|\phi(y)-\phi(x)|>1$,
and that the edge may both be excited or not excited whenever
$|\phi(y)-\phi(x)|=1$.

Consider thus the measure $\gamma_\Lambda(\cdot,\psi)$,
and condition the edge $xy$ to be excited;
the goal is to break symmetry for this excited edge.
Focus on the weight on the left in~\eqref{eq:weight_decomp}.
It is immediate from the definition of $V^*$ that the edge transition distribution $e^{-V^*}\lambda$ can be written 
\[
	e^{-V^*}\lambda=\textstyle\frac12\delta_{-1}+\delta_0+\frac12\delta_1
\]
(where $\delta$ denotes the Dirac measure)
and therefore this distribution equals (up to
normalisation) the distribution of the sum of two independent fair coin flips each valued $\pm1/2$.
If we write $\lambda_*$ for the counting measure on the half-integers,
then 
\[
	e^{-V^*(\phi(y)-\phi(x))}
	=
	\frac12\int 1_{|\phi(x)-z|=\frac12} 1_{|\phi(y)-z|=\frac12}d\lambda_*(z).
\]
 Let $V_*:\mathbb Z+{\textstyle\frac12}\to\mathbb R\cup\{\infty\}$ denote the convex symmetric potential function defined by
 $V_*(x):=\infty\cdot1_{|x|\neq 1/2}$.
 We shall interpret the value of $z$ as a new height associated with the edge 
 $e:=xy$.
 Thus, if we write $\phi(e)$ for $z$, then the previous equation becomes
\[
	e^{-V^*(\phi(y)-\phi(x))}=
	\frac12\int e^{-V_*(\phi(e)-\phi(x))}e^{-V_*(\phi(e)-\phi(y))} d\lambda_*(\phi(e)).
\]
This suggests that we may place a new vertex on the midpoint of $e$---labelled 
also $e$ for convenience---and replace the interaction $V^*$ over the edge $e$
with the interaction $V_*$ over the edges $\{x,e\}$ and $\{y,e\}$.
The marginal distribution of the heights on the original vertex set $\mathbb V$ is invariant under the introduction of this extra height $\phi(e)$.

If an edge is excited, then we may thus associate to it an extra height variable 
which takes values in the half-integers.
In particular, this height is nonzero.
We will exploit this feature by integrating this extra height variable
into an exploration process that leads to the delocalisation result.
Note that conditional on $\phi(x)=\phi(y)=0$, 
the distribution of $\phi(e)$ is that of a fair coin flip
with outcomes $\pm1/2$, independent of all other randomness.
This is the desired symmetry breaking.

Formally, the specification of the further enriched model can be described as follows.
The samples of the model are elements
\[
	(\phi|_\mathbb V,\phi|_\mathbb E,\epsilon)\in\mathbb Z^\mathbb V\times((\mathbb Z+{\textstyle\frac12})\cup\{\infty\})^\mathbb E\times\{0,1\}^\mathbb E,
\]
and the relative weight of each configuration is given by multiplying
\begin{multline*}
	{\textstyle\frac12}\cdot1_{\epsilon(xy)=1}e^{-V_*(\phi(xy)-\phi(x))-V_*(\phi(xy)-\phi(y))}\\
	+1_{\epsilon(xy)=0}1_{\phi(xy)=\infty}\left[e^{-V(\phi(y)-\phi(x))}-e^{-V^*(\phi(y)-\phi(x))}\right]
\end{multline*}
over the edges $xy$ involved in the probability kernel (where $V_*(\infty):=\infty$).
Here $\phi(xy)$ is the intermediate height of the edge,
and the value $\phi(xy)=\infty$ indicates that the intermediate height is not defined.

In the exploration process defined in Subsection~\ref{subsec:finalbigproof},
we will gradually reveal the values of the height function $\phi$ which is
sampled from the measure $\gamma_\Lambda(\cdot,\psi)$.
The measure $\gamma_\Lambda(\cdot,\psi)$---conditioned on the exploration---changes as we reveal more heights and edge states.
As long as factors of the form 
$e^{-V(\phi(y)-\phi(x))}$
appear in the decomposition of this conditioned measure,
we are able to apply the machinery
introduced above: the ideas do not rely on the fact that the sample 
$\phi$ came from the unconditioned measure $\gamma_\Lambda(\cdot,\psi)$.

\subsection{Results from~\cite{SHEFFIELD}}
\label{subsec:quote_sheffield}

For fixed $\Lambda\subset\mathbb V$, write $\mathcal F_\Lambda$
for the $\sigma$-algebra generated by the functions $\phi\mapsto\phi(x)$
with $x$ ranging over $\Lambda$.
Write $\mathcal T$ for the intersection of $\mathcal F_\Lambda$
over all cofinite subsets $\Lambda$ of $\mathbb V$.
Events in $\mathcal T$ are called \emph{tail events}.
A measure is called \emph{extremal} if it satisfies a zero-one law on $\mathcal T$.
We now quote a deep result from the work \emph{Random Surfaces} of Sheffield.

\begin{theorem}[Sheffield,~\texorpdfstring{\cite[Theorem~9.1.1]{SHEFFIELD}}{}]
	Any measure $\mu\in\mathcal G_\mathcal L$ which is ergodic
	is also extremal.
\end{theorem}

The theorem applies in greater generality than is required in this work;
it holds true on any shift-invariant planar graph and for any convex potential function $V$.
The result is phrased slightly differently in~\cite{SHEFFIELD}.
To derive the formulation above, note that
any measure $\mu\in\mathcal G_\mathcal L$ must have the unique gradient Gibbs measure
(with slope zero)
as its gradient,
and ergodicity of $\mu$ implies that the way to choose the heights from the gradients 
does not break extremality.

We shall not employ this result on extremality directly.
Instead, we shall quote two corollaries,
and explain broadly how they are derived.
In the remainder of this subsection, we work in the context of a convex symmetric potential function $V$.
The same corollaries hold true for parity potentials (which are convex \emph{over the odd integers}) after a technical modification which is explained in Subsection~\ref{subsec:proof_main_theorem}.

Fix, throughout this section, some vertex $r\in\mathbb V$,
the \emph{root},
and write $\Lambda_n$ for the set of vertices in $\mathbb V$
which are at a graph distance at most $n$ from $r$.
By the \emph{topology of local convergence} we mean 
the coarsest topology on $\mathcal P(\Omega,\mathcal F)$
that makes the map $\mu\mapsto\mu(A)$ continuous for any 
finite $\Lambda\subset\mathbb V$ and for any $A\in\mathcal F_{\Lambda}$.
It follows from extremality of $\mu$,
that
\[
	\lim_{n\to\infty}
	\gamma_{\Lambda_n}(\cdot,\phi)
	=
	\mu
\]
in the topology of local convergence
for $\mu$-almost every $\phi$ (see~\cite[Lemma~3.2.3]{SHEFFIELD},
or~\cite[Theorem~7.12]{G11} for the same result in the context of general spin systems).

It can be shown that
for any finite set $\Lambda\subset\mathbb V$
and for any admissible height function $\psi$,
the measure $\gamma_\Lambda(\cdot,\psi)$
satisfies the following two properties:
\begin{enumerate}
	\item The distribution of $\phi(r)$ is log-concave,
	\item The measure satisfies the FKG inequality (introduced formally below).
\end{enumerate}
These properties follow directly from the fact that $V$ is a convex potential function,
see~\cite[Lemma~8.2.4]{SHEFFIELD}
and~\cite[Section~9.2]{SHEFFIELD} for details.
The properties are preserved under taking
limits in the topology of local convergence,
which implies the following two results.

\begin{corollary}[Sheffield,~\texorpdfstring{\cite[Lemma~8.2.5]{SHEFFIELD}}{}]
	If $\mu$ is an extremal Gibbs measure,
	then the density of $\phi(r)$ is log-concave.
\end{corollary}

\begin{corollary}[Sheffield,~\texorpdfstring{\cite[Lemma~9.2.1]{SHEFFIELD}}{}]
\label{cor:Sheffield_FKG}
	If $\mu$ is an extremal Gibbs measure,
	then it satisfies the \emph{FKG inequality},
	in the sense that $\mu(fg)\geq\mu(f)\mu(g)$
	for any pair of measurable functions $f,g:\Omega\to[0,1]$
	which are \emph{increasing},
	meaning that $f(\phi)\leq f(\psi)$
	and $g(\phi)\leq g(\psi)$
	for any pair of height functions $(\phi,\psi)$ with $\phi\leq\psi$.
\end{corollary}

Log-concavity of the density of $\phi(r)$ implies in particular
that $\phi(r)$ is integrable, which is exactly the statement 
that we shall aim to contradict.
The FKG property plays a role in the construction of the contradiction.

\subsection{Proof of delocalisation for parity potentials}
\label{subsec:proof_main_theorem}

In this subsection, $V$ is a parity potential,
and $\mathbb G$ is required to be bipartite.
Write $\{\mathbb V_0,\mathbb V_1\}$ for the bipartition
of the vertex set of $\mathbb G$.
We choose the labels such that $\phi$
takes even values on $\mathbb V_0$ and odd values 
on $\mathbb V_1$ almost surely---this is possible
due to ergodicity of $\mu$.
Without loss of generality, we suppose that $r\in\mathbb V_0$.
The two corollaries in the previous subsection remain true,
except that the density of $\phi(r)$ is log-concave \emph{over the even integers}.
In particular, it remains true that $\phi(r)$ is integrable.
The contradiction---and therefore Theorem~\ref{thm:delocalisation_main}---shall follow directly from integrability 
of $\phi(r)$, together with the following lemma.
The remainder of this subsection is dedicated to its proof.

\begin{lemma}\label{lemma:expectation_not_defined}
	For any integer $k$, we have $\mu(\phi(r))\not\in(k-1,k+1)$.
\end{lemma}

For any graph $G$ and some random subset $A$ of its vertices,
write $X_G(A)$ for the event that the complement of $A$ does not contain an infinite cluster.
Consider, for example,
critical site percolation on the triangular lattice,
and write $O$ and $C$ for the set of open and closed 
vertices respectively.
Then both $X_G(O)$ and $X_G(C)$ occur almost surely,
even though neither $O$, nor $C$, percolates.
In the setting of this article, we shall study the events $X_\mathbb G(\{\phi\geq k+1\})$
and $X_\mathbb G(\{\phi\leq k-1\})$.
These events are shift-invariant, and therefore they have probability zero or one for $\mu$ by ergodicity.
They are related to the previous lemma
through the following result.
Note that the lemma holds true also 
when $V$ is any convex symmetric potential,
$\mathbb G$ any locally finite graph,
and $\mu$ any Gibbs measure---as can be seen from the proof.
(The measure $\gamma_\Lambda(\cdot,\psi)$ satisfies the FKG inequality in this more general setting as well.)

\begin{lemma}\label{lemma:exploration}
	Fix $a\in\mathbb Z$ arbitrary.
	If the event $X_\mathbb G(\{\phi\geq a\})$ occurs almost surely,
	then $\mu(\phi(r))\geq a$,
	and if $X_\mathbb G(\{\phi\leq a\})$ occurs almost surely,
	then $\mu(\phi(r))\leq a$.
\end{lemma}

\begin{proof}	
	We prove the first implication for $a=0$;
	the other implications follow by symmetry
	(the parity of $r$ and $a$ does not play a role in the proof).
	The proof is straightforward, and relies on the Markov property 
	and on symmetry and convexity of the potential function $V$.

	For any finite $\Lambda\subset\mathbb V$,
	write $\partial\Lambda$ for the vertices adjacent to $\Lambda$.
	Claim that the expectation of $\phi(r)$ is nonnegative
	in the measure $\gamma_\Lambda(\cdot,\psi)$,
	whenever $\Lambda$ is a finite subset of $\mathbb V$
	containing $r$,
	and $\psi$ an admissible height function 
	with $\psi|_{\partial\Lambda}\geq 0$.
	Indeed, symmetry of the potential function $V$ implies 
	that
	\[
		\gamma_\Lambda(\phi(r),\psi)=-\gamma_\Lambda(\phi(r),-\psi),
	\]
	and the FKG inequality, the Markov property, and nonnegativity of $\psi$ 
	on $\partial\Lambda$ imply that 
	\[
		\gamma_\Lambda(\phi(r),\psi)\geq\gamma_\Lambda(\phi(r),-\psi).
	\]
	This establishes the claim.

	Recall that $\Lambda_n$ denotes the set of vertices at a graph distance at most $n$ from $r$.
	Explore the values of $\phi$ in the following way.
	First reveal the values of $\phi$ on the complement of $\Lambda_n$.
	Then, at each step, select a vertex that has not been revealed
	and which is adjacent to a revealed vertex on which the value of 
	$\phi$ is negative.
	The exploration process terminates when no such vertex can be found.
	This occurs after at most finitely many steps,
	since $\Lambda_n$ is finite.
	Write $R_n=R_n(\phi)$ for the set of vertices that have not been revealed,
	and $A_n=A_n(\phi)$ for the event that $r$ is contained in $R_n$.
	Now
	\[
		\mu(1_{A_n}\phi(r))
		=
		\int 1_{A_n}(\psi)\gamma_{R_n(\psi)}(\phi(r),\psi)d\mu(\psi)
		\geq 0
	\]
	due to the previous claim---note that the values of 
	$\psi$ on $\partial R_n(\psi)$ are nonnegative.
	But $\mu(A_n)\to 1$ as $n\to\infty$
	since the event $X_\mathbb G(\{\phi\geq 0\})$ occurs almost surely.
\end{proof}

To prove Lemma~\ref{lemma:expectation_not_defined},
it suffices to demonstrate that at least one of 
$X_\mathbb G(\{\phi\geq k+1\})$ and $X_\mathbb G(\{\phi\leq k-1\})$ occurs with positive probability for the measure $\mu$.
Without loss of generality, we shall take $k=0$,
and make a number of geometrical observations.
(Note that the parity of $r$ no longer plays a role.)
Recall that $\phi$ takes odd values on $\mathbb V_1$ almost surely.
Let $\mathbb G_1=(\mathbb V_1,\mathbb E_1)$ denote the odd vertex graph: the unique graph 
which has $\mathbb V_1$ as its vertex set, and such that 
two vertices are neighbours if and only if they are at a graph distance 
two in the original graph $\mathbb G$.
The fact that $\mathbb G$ has maximum degree three
implies directly that $\mathbb G_1$ is also a shift-invariant planar graph;
see Figure~\ref{fig:graphs_expanded}.
Moreover, as $\phi$ takes odd values on $\mathbb V_1$,
we observe that this set of vertices can be written 
as the disjoint union of 
$\{\phi\geq 1\}\cap\mathbb V_1$ and $\{\phi\leq-1\}\cap\mathbb V_1$.
Write $\sigma:\mathbb V_1\to\{-1,1\}$
for the unique map such that $\{\sigma=1\}=\{\phi\geq 1\}\cap\mathbb V_1$.
The value of $\sigma(x)$ is called the \emph{spin} of the vertex $x$.

Since each $\mathbb Z$-indexed path through $\mathbb G$ is also a path through 
$\mathbb G_1$ by restricting to half its vertices,
it is immediate that the event
$X_{\mathbb G_1}(\{\sigma=1\})$ is included in the event $ X_\mathbb G(\{\phi\geq 1\})$
and that $X_{\mathbb G_1}(\{\sigma=-1\})$ is included in $X_\mathbb G(\{\phi\leq -1\})$.
Therefore it suffices to demonstrate that 
at least one of $X_{\mathbb G_1}(\{\sigma=1\})$
and $X_{\mathbb G_1}(\{\sigma=-1\})$ has positive probability for $\mu$.
Consider now the converse of this statement.
Then both $\{\sigma=1\}$ and $\{\sigma=-1\}$ must percolate almost surely.
Thus, to arrive at Lemma~\ref{lemma:expectation_not_defined},
it suffices to demonstrate that it is impossible
that both sets percolate simultaneously with positive probability.

For the contradiction we must first derive an intermediate result,
namely that each of $\{\sigma=1\}$ and $\{\sigma=-1\}$ contains at most one infinite cluster almost surely.
This follows from the classical argument of Burton and Keane~\cite{BK}. 

\begin{lemma}
\label{lemma:BK}
The set $\{\sigma=1\}$ contains at most one infinite cluster almost surely.
\end{lemma}

\begin{proof}
	In order to apply the argument of Burton and Keane,
	we must essentially prove insertion tolerance 
	for the percolation $\{\sigma=1\}$.

	Since the law of $\phi$ is shift-invariant 
	and since the distribution of $\phi(x)$ is 
	log-concave for any vertex $x$ (over the odd or even integers,
	depending on the parity of $x$),
	a union bound implies that
	$\|\phi|_{\Lambda_n}\|_\infty/n$
	goes to zero in probability as $n\to\infty$.
	In particular,
	$\liminf_{n\to\infty}\|\phi|_{\Lambda_n}\|_\infty/n=0$
	almost surely.

	Recall that $\partial\Lambda$ denotes the set of vertices adjacent to
	$\Lambda$ in the graph $\mathbb G$.
	For fixed $\phi$ and $n$,
	write $m_n$ for the minimum of $\phi$ on $\partial\Lambda_n$,
	write $\psi_n$ for the function
	\[
		\psi_n:
		\mathbb V\to\mathbb Z\cup\{-\infty\},\,
		x\mapsto\begin{cases}
			-\infty&\text{if $x\not\in\Lambda_{n+1}$,}\\
			m_n+d_\mathbb G(x,\partial\Lambda_n)&\text{if $x\in\Lambda_{n+1}$,}
		\end{cases}
	\]
	and write $\phi_n$ for the height function $\phi\vee\psi_n$.
	It is straightforward to see 
	that 
	\[
		|\phi(y)-\phi(x)|\geq|\phi_n(y)-\phi_n(x)|\in 2\mathbb Z+1
	\]
	for any edge $xy\in\mathbb E$,
	which proves that $\phi_n$ is admissible
	whenever $\phi$ is admissible.
	Moreover, $\phi_n\geq\phi$ by construction.
	The previous paragraph implies furthermore that for fixed $m$,
	the function $\phi_n$ takes strictly positive values on $\Lambda_m$ 
	for arbitrarily large values of $n$ almost surely.
	This construction is similar to a construction in~\cite[Section~14]{LAMMERS}.

	Write $N$ for the number of infinite clusters of $\{\sigma=1\}$.
	Remark that $N$ is deterministic due to ergodicity of $\mu$.
	We aim to demonstrate that $N\leq 1$ almost surely.
	Consider first the case that $2\leq N<\infty$.
	Then for $m$ sufficiently large, the set $\Lambda_m$ intersects all infinite clusters of $\{\sigma=1\}$ with positive probability.
	But, on this event,
	there exists some $n$ such that $\phi_n$ takes strictly positive 
	values on $\Lambda_m$ almost surely.
	This implies that $N=1$ with positive probability
	(by first sampling the height function,
	then resampling its values on $\Lambda_n$),
	a contradiction.
	Finally, consider the case that $N=\infty$ almost surely.
	Then for $m$ sufficiently large, $\Lambda_m$ intersects at least three infinite clusters of $\{\sigma=1\}$ with positive probability.
	On this event, $\phi_n$ takes strictly positive values on $\Lambda_m$ for $n$ sufficiently large almost surely.
	Thus, by reasoning as before, with positive probability $\{\sigma=1\}$ contains 
	an infinite cluster which decomposes into (at least) three infinite clusters 
	when only finitely many vertices are removed (namely those in $\Lambda_n$).
	In other words, with positive probability the set $\Lambda_n$ is a
	\emph{trifurcation box} for $\{\sigma=1\}$ in the sense of the article of Burton and Keane~\cite{BK}.
	Burton and Keane argued that this is impossible,
	because it means that trifurcation boxes would occur with a positive density 
	in the shift-invariant measure $\mu$,
	which in turn leads to a geometrical contradiction.
\end{proof}

Let us collect the intermediate results obtained so far:
\begin{enumerate}
	\item $\mathbb G_1=(\mathbb V_1,\mathbb E_1)$ is a shift-invariant planar graph,
	\item $\sigma$ is an ergodic distribution of spins which satisfies the FKG inequality,
	\item $\{\sigma=1\}$ and $\{\sigma=-1\}$ contain at most one infinite cluster almost surely.
\end{enumerate}
In particular, the FKG inequality in the second statement follows from the fact that each spin of $\sigma$ is an increasing function of $\phi$,
combined with Corollary~\ref{cor:Sheffield_FKG}.
It is known that these intermediate results 
jointly rule out that both $\{\sigma=1\}$
and $\{\sigma=-1\}$ contain an infinite cluster almost surely,
due to Sheffield~\cite[Chapter~9]{SHEFFIELD}.
We also refer to~\cite[Theorem~1.5]{DCRT} for an alternative proof.
Thus, one of the two sets does almost surely not percolate,
say $\{\sigma=-1\}$, in which case 
the event $X_{\mathbb G_1}(\{\sigma=1\})$ occurs.
This establishes the proof of Lemma~\ref{lemma:expectation_not_defined}
(through application of Lemma~\ref{lemma:exploration}),
and consequently that of delocalisation for parity potentials
(Theorem~\ref{thm:delocalisation_main}).

\subsection{Proof of delocalisation for excited potentials}
\label{subsec:finalbigproof}
The letter $V$ denotes an excited potential in this subsection
and $\mathbb G$ is no longer required to be bipartite;
all other notations remain the same.

\begin{lemma}
	For any integer $k$, we have $\mu(\phi(r))\not\in(k,k+1)$.
\end{lemma}

\begin{proof}
	By the same reasoning as in the previous subsection,
	at least one of $X_\mathbb G(\{\phi\leq k\})$
	and $X_\mathbb G(\{\phi\geq k+1\})$ must occur almost surely.
	Lemma~\ref{lemma:exploration} now yields the result. 
\end{proof}

The remainder of this subsection is dedicated to the proof of the following lemma,
which obviously implies the desired contradiction.

\begin{lemma}
	\label{lemma:elaborate_exploration}
	For any integer $k$, we have $\mu(\phi(r))\not\in(k-\frac12,k+\frac12)$.
\end{lemma}

\begin{figure}
    \centering

    \begin{tikzpicture}
        \input{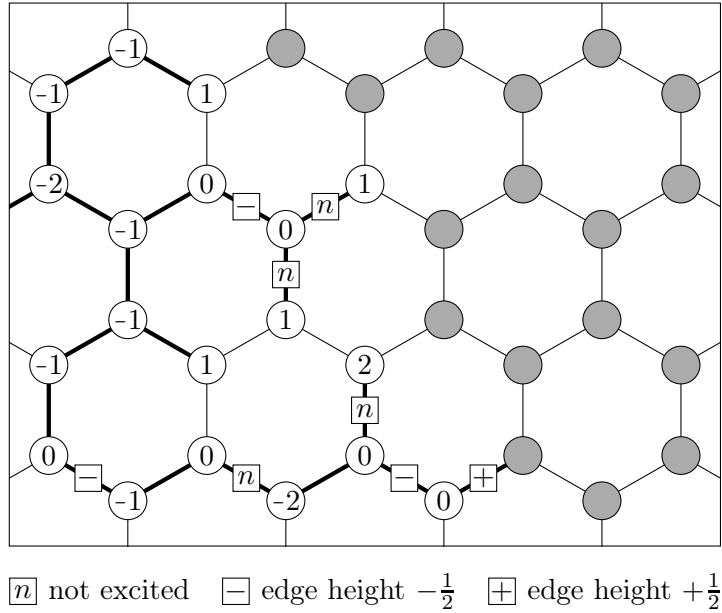}
    \end{tikzpicture}

    \caption{The stopped exploration process, started from the left.
    The thick edges mark the edges that have been considered in the process.
    The heights of the white vertices have been explored;
    the heights of the grey vertices remain unknown.}
    \label{fig:exploration}

    \end{figure}

Set $k=0$ without loss of generality.
Let $\phi$ denote a sample from the Gibbs measure $\mu$,
and let $c:\mathbb E\to\{-1/2,1/2\}$ denote an independent family of fair coin flips.
Write $\sigma:\mathbb E\to\{-1,1\}$ for the unique family 
of spins such that $\sigma(xy)=1$ if and only if $\phi(x)+\phi(y)+c(xy)>0$.
Each spin $\sigma(xy)$ is an increasing function of $(\phi,c)$,
and therefore this family satisfies the FKG inequality.
Remark that the spins are now defined on edges, rather than on vertices as in the previous subsection.
Note that by the same arguments as in the proof of Lemma~\ref{lemma:BK},
the sets $\{\sigma=1\}$ and $\{\sigma=-1\}$ each contain at most one infinite cluster
almost surely.
Since $\mathbb G$ is a shift-invariant planar graph of maximum degree three, we rule out the case that both $\{\sigma=1\}$
and $\{\sigma=-1\}$ contain an infinite cluster almost surely by reasoning as in the previous subsection.
(The phase coexistence result applies because the line graph of $\mathbb G$ is planar;
see Figure~\ref{fig:graphs_expanded}.)
In the following argument we assume that the latter cluster does not percolate
in order to prove that $\mu(\phi(r))\geq 1/2$;
the other case is identical.

We now describe an exploration process which reveals if certain edges are excited or not,
as well as the heights of certain vertices and of all edges which have been revealed to be excited.
Similarly to the proof of Lemma~\ref{lemma:exploration},
we wish the exploration process to terminate before revealing the height of the distinguished root vertex $r$
with high probability.
This is the only place where we require an explicit coupling of the enriched model
with the family of coin flips $c$ introduced above.
We couple the two distributions in the following way:
they are first sampled independently,
and then we set the value of $\phi(xy)$ to $c(xy)$
for each edge $xy$ such that $\phi(x)=\phi(y)=0$.
Clearly this does not change the marginal of the enriched system,
since the distribution of $\phi(xy)$ (conditional on $\phi(x)=\phi(y)=0$) was already that
of an independent fair $\pm\frac12$-valued coin flip.

Recall that $\Lambda_n$ denotes the set of vertices at a distance at most $n$ from the distinguished root vertex $r$.
Sample a height function $\phi$ from $\mu$,
and explore first the values of $\phi$ on the complement of $\Lambda_n$
where $n$ is some large natural number.
Conditional on these values, the law of $\phi$ is given by $\gamma_{\Lambda_n}(\cdot,\phi)$.
Run the following exploration process, which is illustrated by Figure~\ref{fig:exploration}.
Select an edge $e:=xy$ which has not been selected before,
for which the height of $x$ has been explored and $\phi(x)\leq 0$,
and for which the height of $y$ has not been explored.
Immediately explore the value of $\phi(y)$ whenever $\phi(x)<0$.
If $\phi(x)=0$, then we first reveal the state of $e$.
If this edge is not excited, then we immediately reveal the value of $\phi(y)$.
If this edge is excited, then we first explore the intermediate height $\phi(e)$.
If $\phi(e)=-1/2$ then we explore the height of $y$,
and otherwise we choose to not explore the value of $y$ in this step.
Repeat this process until no eligible edges are left.
The process terminates after finitely many steps,
since $\mathbb E(\Lambda_n)$ is finite.

Write $R_n$ for the set of vertices in $\mathbb V$ that have not been revealed
(the grey vertices in Figure~\ref{fig:exploration}).
Let $(x,y)$ denote a (directed) edge in the boundary of $R_n$, in the sense that the value of $\phi(x)$
has been revealed, and the value of $\phi(y)$ not.
Then exactly one of the following two must hold true:
\begin{enumerate}
	\item $\phi(x)>0$ and the state of the edge $xy$ was not revealed,
	\item $\phi(x)=0$, the edge $xy$ was revealed to be excited, and $\phi(xy)=+\frac12$.
\end{enumerate}
Write $E_n^*$ for the set of directed edges of the second type,
and $E_n$ for their undirected counterparts.
Then the law of $\phi$ on $\mathbb V$, conditional on the exploration, is given by
\begin{equation}
\label{eq:condmeasure}
	\frac1Ze^{-H_{R_n,E_n^*}}(\delta_{\phi|_{\mathbb V\smallsetminus R_n}}\times\lambda^{R_n}),
\end{equation}
where $H_{R_n,E_n^*}$ is the Hamiltonian 
\[
	H_{R_n,E_n^*}:=
	\sum_{(x,y)\in E_n^*}V_*(\phi(y)-{\textstyle\frac12})
	+
	\sum_{xy\in\mathbb E(R_n)\smallsetminus E_n}
	V(\phi(y)-\phi(x)).
\]
Note that if an edge $(x,y)$ in the edge boundary of $R_n$ contributes to the sum 
on the right, then it must be of the first type, which automatically means 
that $\phi(x)>0$.
Thus, conditional on this exploration process,
the behaviour of $\phi$ within $R_n$ is that of a random surface
with convex symmetric potential functions
and fixed boundary conditions of at least $1/2$.
Recall from the proof of Lemma~\ref{lemma:exploration}
that the Markov property, the FKG inequality,
and symmetry of the convex potentials $V_*$ and $V$ implies 
that the expectation of $\phi(x)$ in the measure in~\eqref{eq:condmeasure} is at least $1/2$
for each vertex $x\in R_n$.
Conclude that the expectation of $\phi(r)$ is at least $1/2$,
conditional on $r\in R_n$.

Since $\phi(r)$ is integrable for $\mu$,
it suffices to demonstrate that $r\in R_n$ with high probability as $n\to\infty$.
If $r\not\in R_n$,
then some sequence of edges led to the exploration of $r$.
In other words, there must exist a path $(x_k)_{0\leq k\leq m}$ through $\mathbb G$
which starts in $\partial\Lambda_n$
and ends at a neighbour of $r$,
such that $\phi(x_k)\leq 0$ for every vertex,
and such that $\phi(x_kx_{k+1})=c(x_kx_{k+1})=-1/2$
for every edge of this path for which 
$\phi(x_k)=\phi(x_{k+1})=0$.
In particular, this implies that $\sigma(x_kx_{k+1})=-1$ for every edge of this path.
Conclude that $r\in R_n$ whenever 
$r$ is not connected to $\partial\Lambda_n$
for the subgraph $\{\sigma=-1\}\cup\mathbb E(\{r\})$ of $\mathbb G$.
Since $\{\sigma=-1\}$ does not percolate almost surely,
the event $r\in R_n$ has high probability as $n\to\infty$.
This proves that $\mu(\phi(r))\geq 1/2$ whenever the random subgraph $\{\sigma=-1\}$ does not percolate.
The argument for the case that the set $\{\sigma=1\}$ does not percolate is identical.
This establishes the proofs of Lemma~\ref{lemma:elaborate_exploration}
and Theorem~\ref{theorem:delocalisation_even}.

\section*{Acknowledgements}
The author thanks Hugo Duminil-Copin, Matan Harel, Alex Karrila, S\'ebastien Ott, and Martin Tassy for helpful comments and discussions.
The author thanks the referees for their very useful feedback on the first version of the manuscript.

The author was supported by the ERC grant CriBLaM.

\bibliographystyle{amsplain}
\bibliography{references.bib}

\end{document}